\documentclass[12pt]{article}
\usepackage{graphics,pstricks,mfpic}
\usepackage{graphicx}
\usepackage{amsfonts}
\usepackage{amsmath}
\usepackage[all,graph,frame,arrow]{xy}

\title{The Area of a Polygon with an Inscribed Circle}
\author{Marshall W. Buck\thanks{Center for Communications Research, 805 Bunn Dr., Princeton, NJ 08540, buck@idaccr.org} \and
Robert L. Siddon\thanks{2706 Post Oak Ct., Annapolis, MD 21401, rsiddon@aahs.org}}
\date{April 2005}

\newtheorem{proposition}{Proposition}
\begin{document}
\maketitle
\begin{abstract}
Heron's formula states that the area $K$ of a triangle with sides
$a$, $b$, and $c$ is given by
$$
K=\sqrt {s(s-a) (s-b) (s-c)}
$$
where $s$ is the semiperimeter $(a+b+c)/2$.  Brahmagupta, Robbins,
Roskies, and Maley generalized this formula for polygons of up to eight
sides inscribed in a circle. In this paper we derive formulas giving
the areas of any $n$-gon, with odd $n$, in terms of the ordered list
of side lengths, if the $n$-gon is circumscribed about a circle (instead of
being inscribed in a circle). Unlike the cyclic polygon problem,
where the order of the sides does not matter, for the inscribed circle
problem (our case) it does matter. The solution is much easier than
for the cyclic polygon problem, but it does generalize easily to all
odd $n$. We also provide necessary and sufficient conditions for there
to be solutions in the case of even $n$.
\end{abstract}
\newpage
\section{Introduction}\label{sec:1}

Robbins wrote in 1992~\cite{robbinscyclic}:

{\small
\begin{quotation}
Since a triangle is determined by the lengths, $a,b,c$ of its three
sides, the area $K$ of the triangle is determined by these three
lengths.  The well-known formula
\begin{equation*}\label{eq:1.1}
K=\sqrt{s(s-a)(s-b)(s-c)}, \tag{1.1}
\end{equation*}
where $s$ is the semiperimeter $(a+b+c)/2$, makes this dependence
explicit.  (This formula is usually ascribed to Heron of Alexandria,
c.\ 60 BC, although some attribute it to Archimedes.)

For polygons of more than three sides, the lengths of the sides do not
determine the polygon or its area.  However, if we impose the
condition that the polygon be convex and {\it cyclic,} (i.e.,
inscribed in a circle) then the area of the polygon is uniquely
determined.  Moreover, it is a symmetric function of the sides. ...
Given positive real numbers
$a_1,\dots,a_n$, one can construct a convex $n$-gon with the $a_j$'s
as the lengths of the sides provided that the largest $a_j$ is smaller
than the sum of the remaining ones.  In this case it is also possible
to construct a convex cyclic $n$-gon with the same sides (this is not
quite so easy to establish as one might imagine) and this cyclic
$n$-gon has the largest area of all $n$-gons with the given side
lengths. 
\end{quotation}}

Heron's formula has been generalized for cyclic polygons of up to eight
sides \cite{robbinscyclic, maleyrobbinsroskies}.

In this paper, we are interested in finding areas of polygons with $n$ given side
lengths subject to a different constraint. Instead of asking that the
\emph{vertices} of the polygon lie on a circle, we require that the
\emph{sides} of the polygon lie on (are tangent to) a
circle. 

There are several differences with the cyclic polygon problem. First,
in the cases where $n$ is even, there can be infinitely
many solutions. Second, in the case where $n$ is odd,
there are only a finite number of solutions, but different orderings
of the side lengths generally produce solutions with different areas.
In both cases, for there to be any solutions at all, certain equalities and
inequalities in the side lengths $a_j$ must be satisfied.

\section{Heron's Formula}\label{sec:heron}

Let us review the simplest case, $n=3$, producing a simple proof of
Heron's formula that will generalize to $n>3$.

We are given three side lengths $a, b, c$, so the perimeter is
$a+b+c$ and \emph{semi-perimeter} is $s=(a+b+c)/2$.
For there to be a triangle with those side lengths, the following
(triangle) inequalities must be satisfied:
\begin{align*}
a &\leq b+c \\
b &\leq c+a \\
c &\leq a+b.
\end{align*}
For there to be an inscribed circle, and for the area of the triangle
to be nonzero, we need strict inequalities above.

\begin{figure}[h!]
\begin{center}
\includegraphics[scale=0.7]{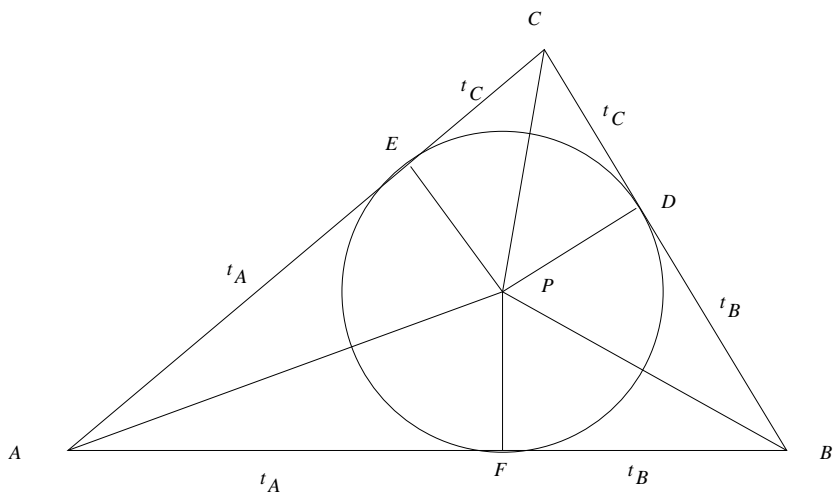}
\end{center}
\caption{Diagram for proof of Heron's formula}
\label{fig:n3}
\end{figure}

Refer to Figure~\ref{fig:n3}.
Denote by $A$ the vertex opposite the side of length $a$, by $B$ the
vertex opposite the side of length $b$, and by $C$ the vertex opposite
the side of length $c$. Thus, the side $BC$ has length $a$,
(more succinctly, $|BC|=a$); furthermore, $|AC| =b$ and $|AB|=c$.
Let $P$ be the center of the inscribed circle, and let $D,E,F$ be the
feet of perpendiculars dropped from $P$ to the sides $BC,CA,AB$
respectively. (Thus, $D,E,F$ are the points of tangency of the
triangle with the inscribed circle.)
These perpendiculars are radii of the circle, and all have the same
length $r$. The area $K$ of the triangle $\triangle ABC$ is the sum of the areas
of the three triangles $\triangle ABP$, $\triangle BCP$, and
$\triangle CAP$, having respective
areas $cr/2$, $ar/2$, and $br/2$, summing to $(a+b+c)r/2$. Thus, we
have the well-known formula relating the radius of the inscribed circle,
the area of the triangle, and the semi-perimeter $s$:
\begin{equation}\label{eq:K=rs}
K = s r.
\end{equation}

The two tangents from $A$ to the circle, $AE$ and $AF$, have the same
length $t_A$. Similarly, the tangents from $B$ have length $t_B$, and
the tangents from $C$ have length $t_C$.
Meanwhile, each side of the triangle is composed of two of these
tangent lines.
Thus, we have
\begin{eqnarray*}
a &=& |BC| = |BD| + |DC| = t_B + t_C \\
b &=& |AC| = |AE| + |EC| = t_A + t_C \\
c &=& |AB| = |AF| + |FB| = t_A + t_B \\
a+b+c &=& 2(t_A+t_B+t_C)\\
s &=& t_A + t_B + t_C \\
t_A &=& s - a\\
t_B &=& s - b\\
t_C &=& s - c.
\end{eqnarray*}

Let $\theta_A, \theta_B, \theta_C$ denote the angles $\angle
APE=\angle APF$, $\angle BPD=\angle BPF$, and $\angle CPD =\angle CPE$
respectively. We must have $\theta_A + \theta_B + \theta_C=\pi$.
Meanwhile, we observe that $\theta_A, \theta_B, \theta_C$ are the
arguments of the three complex numbers $r+it_A, r+it_B, r+it_C$, respectively.
The argument of the complex number $(r+it_A)( r+it_B)( r+it_C)$ must
then be $\theta_A + \theta_B + \theta_C=\pi$, which means that that
the imaginary part of $(r+it_A)( r+it_B)( r+it_C)$ must be zero.
But
\[
\Im\left\{(r+it_A)( r+it_B)( r+it_C)\right\}= (t_A+t_B+t_C) r^2 - t_A
t_B t_C, 
\]
so
\begin{equation}\label{eq:r}
r^2 = \frac{t_A t_B t_C}{t_A+t_B+t_C} = \frac{(s-a)(s-b)(s-c)}{s}.
\end{equation}
Combining \eqref{eq:K=rs} and \eqref{eq:r}, we have Heron's formula
\[
K^2 = s^2 r^2 = s^2 \frac{(s-a)(s-b)(s-c)}{s} = s(s-a)(s-b)(s-c).
\]

\section{Generalization to \protect{\boldmath{$n>3$}}: Finding the Tangent Lengths}

Our solution involves two steps: (1) determine the tangent lengths
from the side lengths; (2) determine the radius from the tangent
lengths. The first step depends upon the ordering of the side lengths,
but the second step produces a polynomial equation of degree $\lfloor(n-1)/2\rfloor$
satisfied by $r^2$, where the coefficients of the polynomial are
elementary symmetric functions of the tangent lengths.

Denote the side lengths in cyclic order by $a_1,a_2,\ldots,a_n$, and the
tangent lengths in order by $t_1,t_2,\ldots,t_n$, designated so that
we have
\begin{align}\label{eq:u=t}
a_1 &= t_1 +t_2 \\ \nonumber
a_2 &= t_2 +t_3 \\ \nonumber
&\vdots \\ \nonumber
a_n &= t_n +t_1 .
\end{align}
Again, we define the semi-perimeter 
\begin{eqnarray*}
s &=& (a_1+a_2+\cdots+a_n)/2 \\
  &=& t_1+t_2+\cdots+t_n,
\end{eqnarray*}
and \eqref{eq:K=rs} still holds in the general case.

If $n$ is odd, we can determine the tangents $t_i$ from the side
lengths by inverting the linear equations above, getting
each $t_i$ as an alternating sum:
\begin{align*}
t_1 &= (a_1 - a_2 + a_3 - \cdots + a_n)/2\\
t_2 &= (a_2 - a_3 + a_4 - \cdots -a_n + a_1)/2\\ \nonumber
    &\vdots \\ 
t_n &= (a_n - a_1 + a_2 - \cdots + a_{n-1})/2.
\end{align*}
In this case, the condition that there be a polygon circumscribing a
circle is that each of
the $t_i$ is positive, which converts into the following $n$ inequalities
in the side lengths:
\begin{equation}
a_1 + a_3 + a_5 + \cdots + a_n > a_2 + a_4 + \cdots + a_{n-1}
\end{equation}
and all other cyclic rotations of the indices.

If $n$ is even, we wish to find necessary and sufficient conditions on
the $a_i$ values and ordering, for there to be a corresponding set of
positive tangent lengths satisfying the equations \eqref{eq:u=t}.
First, we can write down equations determining each $t_i$ in terms of
a single one, $t_1$:
\begin{eqnarray*}
t_2 &=& a_1 - t_1 \\
t_3 &=& a_2 - t_2 = a_2 -a_1 + t_1 \\
t_4 &=& a_3 - t_3 = a_3 -a_2 + a_1 - t_1 \\
&\vdots& \\
t_n &=& a_{n-1} - t_{n-1} = a_{n-1} - a_{n-2} +\cdots + a_1 - t_1\\
t_1 &=& a_n - a_{n-1} + \cdots + a_2 - a_1 + t_1.
\end{eqnarray*}
For there to be any solution at all, it is necessary that the sum of
the $a_i$ for odd $i$ is equal to the sum of the $a_i$ for even $i$.
Beyond that, it is necessary that we be able to choose $t_1$ to make
all of the $t_i$ positive. 
We can do that if and only if
we can choose $t_1$ less than all of the following:
\begin{align*}
 &a_1 \\
 &a_1 - a_2 + a_3\\
 &a_1 - a_2 + a_3 - a_4 + a_5 \\
 &\vdots \\
 &a_1     - a_2     + \cdots -a_{n-2} + a_{n-1},
\end{align*}
and, simultaneously, greater than all of the following:
\begin{align*}
 &0\\
 &a_1 - a_2 \\
 &a_1 - a_2 + a_3 - a_4  \\
 &\vdots \\
 &a_1     - a_2     + \cdots +a_{n-3} -a_{n-2}.
\end{align*}
In other words, we need
\begin{eqnarray*}
&&\hspace{-.3in}\min \left( a_1,
 a_1 - a_2 + a_3,
 a_1 - a_2 + a_3 - a_4 + a_5,
 \dots,
 a_1     - a_2     + \cdots -a_{n-2} + a_{n-1}\right)> \\
&&\hspace{-.3in}\max \left(0,
 a_1 - a_2 ,
 a_1 - a_2 + a_3 - a_4  ,
 \dots ,
 a_1     - a_2     + \cdots + a_{n-3} - a_{n-2}
\right).
\end{eqnarray*}
This happens if the alternating sum of the $a_i$ of length $n$
vanishes, and if every contiguous, odd length, alternating sum
is positive.  Of course, we are given that the length one sums, like
$a_2$, are positive. But, we also need (if $n\geq 6$) the additional
alternating sum conditions like $ a_1 - a_2 + a_3 >0$.  Not all of
these are needed, even of those length $n/2$ or less. The vanishing of
the alternating sum of length $n$ allows us to derive the alternating
sum condition for a chain of length $k$ from the complementary chain
of length $n-k$. If we count the conditions $a_j>0$, there are a total
of $n^2/4$ independent inequalities needed, and these can be taken to
be alternating sums of any odd length, starting from an $a_j$ for odd
$j$. It is easy to show that all the inequalities are
necessary. Special cases of this result for $n=4,6,8$ are stated below.

\begin{proposition} 
There will be a circle-enclosing
quadrilateral with given side lengths $a_1,a_2,a_3,a_4$ if and only if the
equality
\[
a_1 + a_3 = a_2 + a_4,
\]
and the following four inequalities are satisfied:
\[
a_1>0, a_2>0, a_3>0, a_4>0.
\]
\end{proposition}

\begin{proposition} 
There will be a circle-enclosing
hexagon with given side lengths $a_1,a_2,\ldots,a_6$ if and only if the
equality
\[
a_1 + a_3 + a_5= a_2 + a_4 + a_6,
\]
and the following nine inequalities are satisfied:
\begin{align*}
a_1 &>0, a_2>0, \ldots, a_6>0 \\
a_1 &- a_2 + a_3 > 0\\
a_3 &- a_4 + a_5 > 0\\
a_5 &- a_6 + a_1 > 0.
\end{align*}
\end{proposition}

\begin{proposition} 
There will be a circle-enclosing
octagon with given side lengths $a_1,a_2,\ldots,a_8$ if and only if the
equality
\[
a_1 + a_3 + a_5 + a_7= a_2 + a_4 + a_6 +a_8,
\]
and the following $16$ inequalities are satisfied:
\begin{align*}
a_1 &>0, a_2>0, \ldots, a_8>0 \\
a_1&- a_2 + a_3 > 0\\
a_2&- a_3 + a_4 > 0\\
   &\vdots\\    
a_7&- a_8 + a_1 > 0\\
a_8&- a_1 + a_2 > 0.
\end{align*}
\end{proposition}

Computationally, rather than check the large number of inequalities
($n^2/4$ such) it is better to proceed as follows. From left to
right, incrementally calculate the sum $ v=a_1-a_2+a_3-\cdots-a_n
$. 
We calculate the allowed range for $t_1$ as the interval $(l,r)$,
starting with the interval $(0,a_1)$ and then executing the following
pseudo C-code:

{\small
\begin{verbatim}
   double u[n+1]; // u[1],...,u[n] are side lengths
   double r, ell, v=0.0; /* (ell,r) will be the choice interval for t_1
   ell = 0.0; r = u[1]; /* initialize choice interval to (0,u_1))
   for(int i = 1; i <= n; i++){
      if( i&1 ) { /* i is odd */
         v += u[i];
         if( v < r ) r = v;
      } else { /* i is even */
         v -= u[i];
         if( v > ell ) ell = v;
      }
   }
   if( v != 0.0 ) printf("The alternating sum is not zero as required\n");
   if( ell >= r ) printf("The choice interval for t_1 is empty\n");
\end{verbatim}
}
(In an actual program, we would not compare a floating point number to ``0.0''.)

\section{Generalization to \protect{\boldmath{$n>3$}}: from Tangent Lengths to the Radius of
the Inscribed Circle}

Given the tangent lengths $t_1,t_2,\ldots,t_n$ (summing to the
semiperimeter $s$, and satisfying the edge length constraints
$a_1=t_1+t_2$, $a_2=t_2+t_3$, etc.) we may now construct the polygon
about any circle of radius $r$, by using complex numbers
$q_1=r,q_2,\ldots,q_n$ to represent
the points of tangency and the complex numbers $p_1,p_2,\ldots,p_n$
to represent the vertices of the polygon.
From each $q_j$ we produce the next $p_j$ point, in the counterclockwise
direction, and having tangent length $t_j$, by multiplying $q_j$ by
$(r+it_j)/r$.
Similarly, we calculate the next $q_{j+1}$ by multiplying $p_j$
by $r/(r-it_j)$. Thus, we have
\begin{eqnarray*}
q_1 &=& r\\
p_1 &=& \frac{r+it_1}r q_1 \\
q_2 &=& \frac{r}{r-it_1} p_1\\
p_2 &=& \frac{r+it_2}r q_2\\
q_3 &=& \frac{r}{r-it_2} p_2 \\
p_3 &=& \frac{r+it_3}r q_3\\
&\vdots& \\
p_n &=& \frac{r+it_n}r q_n\\
q_1 &=& \frac{r}{r-it_n} p_n.
\end{eqnarray*}
For the polygon to close, it is necessary and sufficient that
\[
\prod_{j=1}^n \frac{r+ i t_j}{r- i t_j} =1.
\]
This happens if and only if the product
$$
\prod_{j=1}^n (r+ i t_j)
$$
is real.
Define the elementary symmetric functions
$\sigma_1,\sigma_2,\ldots,\sigma_n$ of the tangent lengths
$t_1,t_2,\ldots,t_n$ as usual.
In particular, $\sigma_1 = t_1+t_2+\cdots+t_n = s$ and $\sigma_n = t_1 t_2 \cdots
t_n$.  Also, define $\sigma_0=1$.
These can be calculated from the identity
\[
\prod_{j=1}^n (x + t_j) = \sum_{j=0}^n \sigma_j x^{n-j}.
\]
Rewriting the identity, we obtain
\[
\prod_{j=1}^n (r+ i t_j) = \sum_{j=0}^n i^{j} \sigma_j r^{n-j}.
\]
Meanwhile, we want the imaginary part of the left-hand side to vanish, so we
can rewrite our condition as
\begin{equation}\label{eq:sigr}
\sigma_1 r^{n-1} - \sigma_3 r^{n-3} +\cdots = 0.
\end{equation}
When $n=2k+1$ is odd, $r^2$ must be a root of the polynomial 
\begin{equation}
\sigma_1 x^k - \sigma_3 x^{k-1} + \cdots + (-1)^k \sigma_n.
\end{equation}
When $n=2k+2$ is even, we can take out one factor of $r$ from
\eqref{eq:sigr}, and say that $r^2$ must be a root of the polynomial 
\begin{equation}
\sigma_1 x^{k} - \sigma_3 x^{k-1} + \cdots + (-1)^k \sigma_{n-1}.
\end{equation}
In both cases there will be exactly $k$ roots
$r_1^2>r_2^2>\cdots>r_k^2$. The largest radius, $r_1$, will give rise to a
convex polygon enclosing the circle of radius $r_1$. 
The other roots will have the polygon going around the circle multiple
times. For example, for $n=5$ and $k=2$, the radius $r_2$ will be that
of a circle
inscribed in a pentagonal star.

One way to see that there are $k$ solutions is to consider the sum of
the central angles that needs to be a multiple of $\pi$. Let
$$
f(r)=\tan^{-1}\frac{t_1}r + \tan^{-1}\frac{t_1}r+\cdots+\tan^{-1}\frac{t_n}r.
$$
The continuous function $f(r)$ on the interval $r\in (0,\infty)$ converges to $\frac{n\pi}2$ as $r \rightarrow 0$,
converges to 0 as $r\rightarrow \infty$, and monotonically decreases
as $r$ goes from 0 to $\infty$.
Thus, by the intermediate value theorem, $f$ takes on the values
$\pi,2\pi,\ldots,k\pi$ at $r=r_1,r=r_2,\ldots,r=r_k$.

Since general polynomial equations of degree up to four are solvable in
radicals, we can similarly solve our problem for  $k\leq 4$, thus for $n\leq 10$.

\section{Quadrilaterals}

Given four sides $a_1,a_2,a_3,a_4$, the largest area of a
quadrilateral is for the cyclic quadrilateral, and that area is given
by Brahmagupta's formula
\begin{equation}
K = \sqrt{(s-a_1)(s-a_2)(s-a_3)(s-a_4)}.
\end{equation}

Meanwhile, if $a_1 + a_3 = a_2 + a_4$, we will be able to inscribe
circles in various ways. We would like to find the quadrilateral of
largest area, among those containing a circle tangent to all four
sides.

In this case, because of the constraint on the $a_j$'s, Brahmagupta's
formula simplifies.
For example, 
\[s-a_1 = (-a_1+a_2+a_3+a_4)/2= (-a_1+a_2+a_1+a_2)/2=a_2.\]
So,
\begin{equation}
K = \sqrt{a_1 a_2 a_3 a_4}
\end{equation}
in this special case.
(See also \cite[Exercise 1, page 60]{coxeter-greitzer}.)

If we set the tangent lengths to be
\begin{eqnarray}\label{eq:soln}
t_1 &=& \frac{a_1 a_4}{a_1+a_3}\\ \nonumber
t_2 &=& \frac{a_2 a_1}{a_1+a_3}\\ \nonumber
t_3 &=& \frac{a_3 a_2}{a_1+a_3}\\ \nonumber
t_4 &=& \frac{a_4 a_3}{a_1+a_3}
\end{eqnarray}
it is easy to check that the conditions $a_1=t_1+t_2$, etc. are
satisfied (making use of $a_1+a_3 = a_2+a_4$).  Furthermore, one can verify
that $t_1 t_3 = t_2 t_4$, and $s=a_1+a_3 = a_2+a_4$.
Now we can calculate the area using our earlier results.
First, $r^2$ satisfies $ \sigma_1 r^2 = \sigma_3 $, $\sigma_1 = s$,
and
\begin{eqnarray*}
\sigma_3 &=&  t_1 t_2 t_3 +  t_1 t_2 t_4 +  t_1 t_3 t_4 + t_2 t_3 t_4 \\
 &=&  t_1 t_3 (t_2 + t_4) + t_2 t_4 ( t_1+t_3) \\
 &=&  t_1 t_3 (t_1 + t_2 + t_3 + t_4)\\
 &=&  s t_1 t_3 \\
 &=&  s^{-1} a_1 a_2 a_3 a_4,
\end{eqnarray*}
so we have
\begin{eqnarray*}
K^2 &=& r^2 s^2 = s^2 (\sigma_3/\sigma_1)\\
&=& s \sigma_3 = a_1 a_2 a_3 a_4.
\end{eqnarray*}
This construction of the quadrilateral achieves the maximum area, so
it simultaneously is inscribed in a circle and has a circle inscribed
in it. 


Samuel Kutin~\cite{kutin} remarks that the solution
\eqref{eq:soln} can be derived simply as follows. He notes that in the
case of a quadrilateral inscribed in a circle, opposite internal
angles are supplementary, so opposite half internal angles are
complementary. If there is an inscribed circle of radius $r$, then
there are similar right triangles telling us that $t_1/r = r/t_3$ and
$t_2/r= r/t_4$, and therefore $t_1 t_3 = t_2 t_4$. But, from this
condition alone, we can then derive \eqref{eq:soln}. For example,
\begin{eqnarray*}
a_1 a_2 &=& (t_1+t_2)(t_2+t_3) \\
        &=& t_1 t_2 + t_1 t_3 + t_2^2 + t_2 t_3\\
        &=& t_1 t_2 + t_2 t_4 + t_2^2 + t_2 t_3 \\
        &=& t_2( t_1 + t_4 + t_2 + t_3) \\
        &=& t_2( (t_1 + t_2) + (t_3 + t_4)) \\
        &=& t_2( a_1 + a_3),
\end{eqnarray*}
giving $t_2 = a_1 a_2/(a_1+a_3)$, the second line of \eqref{eq:soln}.

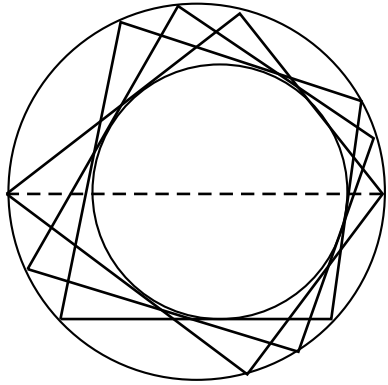
\begin{figure}
\moveright.8in\vbox{
\begin{picture}(100,100)
\put(146.7,-0.75){\pscircle{48.6pt}}
\put(138,-0.75){\pscircle{71.6pt}}
\put(66,-30){\psline[linestyle=dashed,linewidth=1pt]{-}(5.0,1)(0,1)}
\put(66,-30){\psline[linestyle=solid,linewidth=1pt]{-}(0,1)(3.1,3.4)(5,1)(3.2,-1.4)(0,1)}
\put(74,-30){\psline[linestyle=solid,linewidth=1pt]{-}(0,0)(2.0,3.5)(4.6,1.74)(3.6,-1.1)(0,0)}
\put(58,-49){\psline[linestyle=solid,linewidth=1pt]{-}(1,0)(4.6,0)(5.0,2.9)(1.8,3.95)(1,0)}
\end{picture}}
\vglue1in
\caption{Construction of quadrilaterals via Poncelet's Porism}
\label{fig:circle}
\end{figure}

Stephen DiPippo~\cite{dipippo} gives a nice way to construct examples
(see Figure~\ref{fig:circle}). He starts with any right triangle,
makes a mirror image copy of the triangle about the hypotenuse, joins
the two triangles to form a quadrilateral with an inscribed circle.
Meanwhile, the hypotenuse is a diameter of a circle containing all
four vertices of the quadrilateral. This makes one example of such a
quadrilateral, but it is ``too symmetric''. Then he remarks that one
can apply Poncelet's Porism~\cite{griffiths-harris}, and know that one can start at \emph{any}
point on the outer circle, drop a tangent to the inner circle,
intersect that tangent with the outer circle, drop another tangent,
and continue, forming a new quadrilateral having the same pair of
circles on the inside and outside. By such means, he constructs a
general example.

\section*{Keywords}
{\small{\tt polygon, cyclic, circle, inscribed, circumscribed, area, radius, \break}
{\tt Heron, Brahmagupta, Robbins, Maley, Litman, Roskies, Hero}}
\end{document}